\documentclass{amsart}
\usepackage[margin=1.5in]{geometry}

\usepackage[all]{xy}
\usepackage{tikz-cd}
\usepackage{enumerate}
\usepackage{amsfonts}
\usepackage{amssymb, amscd}
\usepackage{graphicx}
\usepackage{mathtools}
\usepackage{xcolor}
\usepackage{float}

\newcommand{\sn}{\smallskip\noindent}

\newtheorem{teo}{Theorem}[section]

\theoremstyle{definition}

\title{A new and simple proof of the false centre theorem}

\author[L. Montejano]{Luis Montejano}
\address{UNAM at Quer\'etaro. Mexico}

\author[E. Morales-Amaya]{Efren Morales-Amaya}
\address{Facultad de matem\'aticas at Acapulco, UAGro. Guerrero, Mexico.}

\begin{document}
\maketitle
\section{Introduction and preliminaries}

We say that a set $A\subset \mathbb R^n$ is \emph{symmetric} if and only if there is a translated copy $A'$ of $A$ such that 
$A'=-A'$. In this case, if  $A'=A-x_0$, we say that $x_0$ is the center of symmetry of $A$. By convention, the empty set $\emptyset$ is symmetric.

The purpose of this paper is to give a new and simple proof of the following theorem:

\begin{teo}\label{Main}[False Centre Theorem]
Let $K$ be a  convex body in euclidean $3$-space and let ${p}$ be a point of $\mathbb R^3$.  Suppose that for every plane $H$ through $p$,
the section $H\cap K$ is symmetric. Then either $p$ is a centre of symmetry of $K$ or $K$ is an ellipsoid.  
\end{teo}

For the proof we will use the following two known theorems. 

\begin{teo}\label{teoss}
Let $K$ be a convex body in euclidean $3$-space.  Suppose that for every plane $H$,
the section $H\cap K$ is symmetric. Then $K$ is an ellipsoid.  
\end{teo}

\begin{teo} \label{teorogers}
Let $K$ be a convex body in euclidean $3$-space and let ${p}$ be point in $\mathbb R^3$.  Suppose that for every plane $H$ through $p$,
the section $H\cap K$ is symmetric. Then $K$ is symmetric.  
\end{teo}

Theorem \ref{Main} was first proved in all its generality by D. G. Larman \cite{L}. Theorem \ref{teorogers} was proved by C. A. Rogers \cite{R}, when $p\in$ int $K$ and by  G. R.  Burton in general (Theorem 2 of \cite{B1}). Theorem \ref{teoss} was first proved  by H. Brunn \cite{Bru} under the hypothesis of regularity and in general by G. R. Burton \cite {B2} (see (3.3) and (3.6) of  Petty's survey \cite{P}). For more about characterization of ellipsoids see \cite{So} and Section 2.12 of \cite{MMO}.  

We need some notation. 
Let $K$ be a  convex body in euclidean $3$-space $\mathbb R^3$, let $p\in \mathbb R^3$ and let $L$ be a directed line. 
We denote by  $p_L$ the directed chord  $(p+L)\cap K$ of $K$, by $|p_L|$ the length of the chord $p_L$ and  by 
$\frac{p}{L}\in \mathbb R \cup \{\infty\}$, the radio in which the point $p$ 
divides the directed interval $p_L$. That is, if $p_L=[a,b]$, then $\frac{p}{L}=\frac{pa}{bp}$, where $pa$ and $bp$ denotes the signed length of the directed chords $[p,a]$ and $[b,p]$ in the directed line $p+L$, and by convention, if $p_L=[p,p]$, then $\frac{p}{L}=1$.  If $p_L=(p+L)\cap K=\emptyset$, then by convention $|p_L|=-1$ and $\frac{p}{L}=-1.$ 

Suppose now $B$ is a convex figure in the plane and  let $p$ be a point of $\mathbb R^2.$
If $B$ is symmetric with center the origin, $L$ is a directed line through the origin and $q=-p$, then:
\begin{enumerate}
\item  $|p_L|= |q_L|$, and
\item  $\frac{q}{L}\frac{p}{L}=1$, 
\end{enumerate}
\noindent where by convention $0\infty=\infty0=1$.

Conversely, if $B$ is a convex figure and $p,q\in \mathbb R^2$ are two  points for which $(1)$ and $(2)$ holds, for every directed line $L$ through the origin,
then $B$ is symmetric with center at the midpoint of $p$ and $q$. Essentially, this is so because, if $p=-q$ and $p_L=[a,b]$, then 
$q_L=[-b,-a]$.

\section{The proof of Theorem \ref{Main}}  

Let $K$ be a convex body in euclidean $3$-space and let ${p}$ be point of $\mathbb R^3$. Suppose that for every plane through $p$,
the corresponding section is symmetric. By Theorem \ref{teorogers}, we may assume that $K$ is symmetric with the center at the origin.  Suppose $p$ is not the origin $0$. We shall prove that $K$ is an ellipsoid.  Let $H$ be a plane through the origin that does not contain the point $p$.  By hypothesis, the plane $(p+H)\cap K$ is symmetric.  Suppose $v\not=p$ is the center of $(p+H)\cap K$, and $w=tv$, for some $t\in \mathbb R$. We shall prove first that 
$(w+H)\cap K$ is symmetric with the center at $w$. 
For that purpose, it will be enough to prove that for every line $L\subset H$ through the origin:
\begin{enumerate}

\item $|(p-v+w)_L|= |(v-p+w)_L|$, and
\item $\frac{v-p+w}{L}\frac{p-v+w}{L }=1.$
\end{enumerate}

Note that the points $p-v+w$ and $v-p+w$ lie in $(w+H)\setminus \partial K$ and if $L\subset H$, then $(p-v+w)_L$ and $(v-p+w)_L$ are chords of $(w+H)\cap K$.
 
 Let $\Gamma$ be the plane through the origin generated by $L$ and $v$. Hence, by hypothesis, $(p+\Gamma)\cap K$  is symmetric. 
Suppose first that $K$ is strictly convex. In order to prove $(1)$ and $(2)$, we shall prove first:

\begin{itemize}
\item[a)]  $|p_L|=|(p-2v)_L|,$
\item[b)]  the center of $(p+\Gamma)\cap K$ lies in $(p-v)_L$, and
\item[c)]  $\frac{p}{L}=\frac{p-2v}{L}$.
\end{itemize}

 Since $(p+H)\cap K$ is a symmetric section with center at $v$, then $|p_L|=|(2v-p)_L|$ and $\frac{p}{L}\frac{2v-p}{L}=1$. By the symmetry of $K$, 
$|(p-2v)_L|=|(2v-p)_L|$ and $\frac{p-2v}{L}\frac{2v-p}{L}=1$. Consequently, $|p_L|=|(p-2v)_L|$. So both chords, $p_L$ and $(p-2v)_L$, of the symmetric section 
$(p+\Gamma)\cap K$ have the same length. By the strictly convexity of $K$, the parallel mid chord contains the center, that is, the center of $(p+\Gamma)\cap K$ lies in $(p-v)_L$. Furthermore, $\frac{p}{L}=\frac{p-2v}{L}$. This proves a), b) and c).

\begin{figure}
\includegraphics[scale=.5]{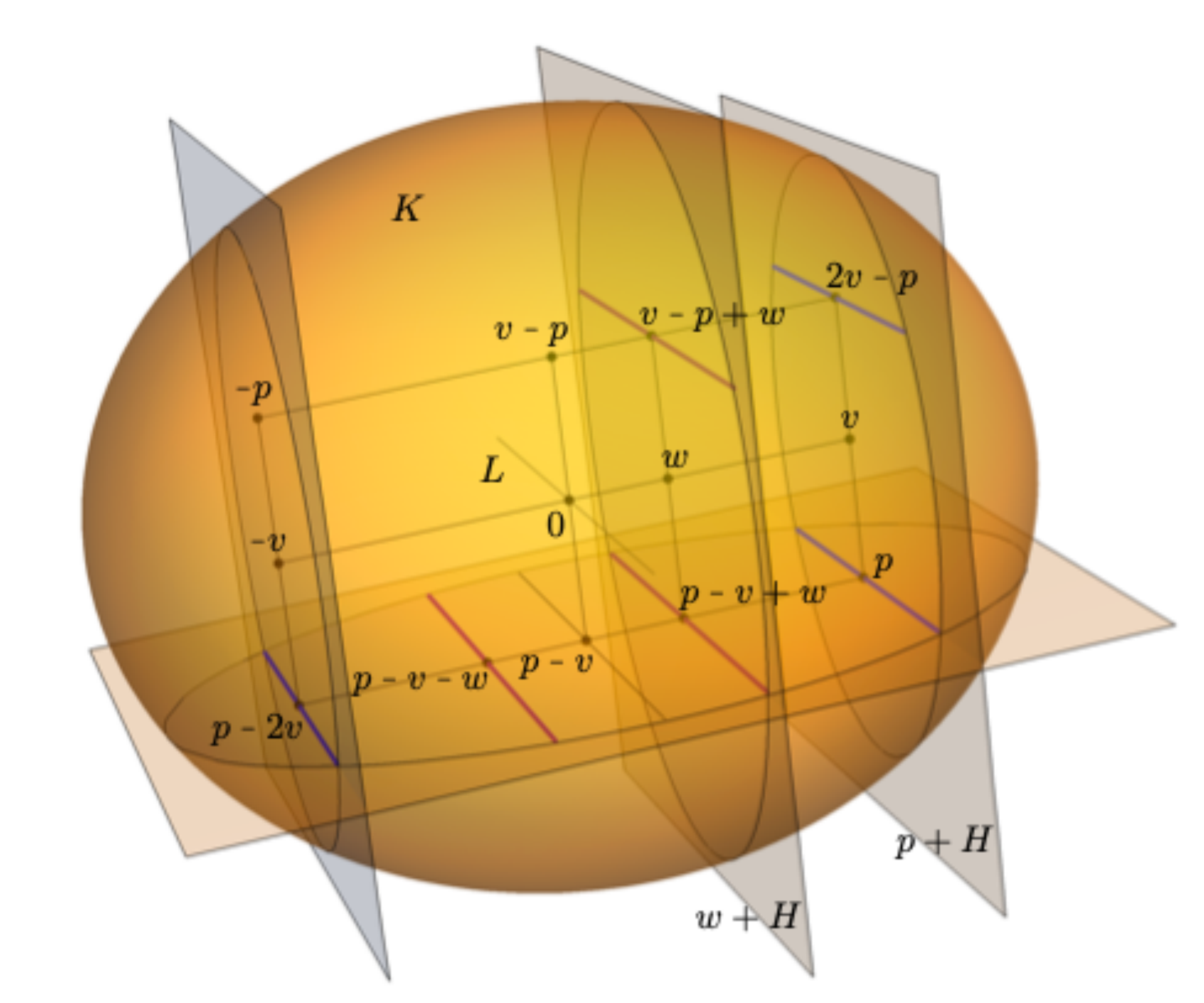}
\caption{}
\end{figure}

Let us prove now that $a), b)$ and $c)$ implies $(1)$ and $(2)$.  The parallel chords $(p-v-w)_L$ and $(p-v+w)_L$ of the symmetric section 
$(p+\Gamma)\cap K$  have as a mid chord $(p-v)_L$. The fact  that the center of $(p+\Gamma)\cap K$ lies in $(p-v)_L$ implies that  $|(p-v+w)_L|= |(p-v-w)_L|$. Since $K$ is symmetric with center at the origin, then (1) holds, that is  $|(p-v+w)_L|= |(v-p+w)_L|$.  Since $\frac{p}{L}=\frac{p-2v}{L}$, then 
$(p-2v)_L+2v=p_L$. By symmetry of $(p+\Gamma)\cap K$, $(p-v-w)_L+2w=(p-v+w)_L$. Hence $\frac{p-v-w}{L}=\frac{p-v+w}{L}$. On the other hand, by the symmetry of $K$, $\frac{p-v-w}{L}\frac{v-p+w}{L}=1$. Consequently, $\frac{v-p+w}{L}\frac{p-v+w}{L}=1$, thus proving $1)$ and $2)$ and hence that $(w+H)\cap K$ is symmetric with the center at $w$.

This proves that every section of $K$ parallel to $H$ is symmetric.  
Let us prove now that the collection of planes through the origin such that the center of the section $(p+H)\cap K$ is not $p$, is dense. 
Let $\Omega$ be the collection of planes through the origin such that the center of the section $(p+H)\cap K$ is $p$ and let $H\in$ int$\Omega$. If this is the case, by symmetry of $K$, the center of the section $(-p+H)\cap K$ is $-p$ and therefore, $\big((-p+H)\cap K\big)+2p = (p+H)\cap K$. Since the same hold for every section sufficiently close, we conclude that  $(H\cap \partial K)+\{tp\in \mathbb R^3\mid|t|\leq 1\}\subset \partial K$, contradicting the strictly convexity assumption. Consequently, the collection of symmetric sections of $K$ is dense and 
since the limit of a sequence of symmetric sections is a symmetric section, then every section of $K$ is symmetric.  By Brunn's Theorem \ref{teoss}, $K$ is an ellipsoid. 
In the non strictly convex case, our arguments for the case in which the center of $(p+H)\cap K$ is not $p$, 
only showed that $(w+H)\cap K$ is symmetric, when $w=tv$ from $|t|<1$, thus proving that every section of $K$ sufficiently close to the origin and parallel to $H$ is symmetric and hence that $H\cap\partial K$ is contained in a shadow boundary.  
On the other side, if $H\in$ int$\Omega$ and 
$(H\cap \partial K)+\{tp\in \mathbb R^3\mid|t|\leq 1\}\subset \partial K$, then clearly $H\cap\partial K$ is contained in a shadow boundary.
Consequently, by Blaschke's Theorem 2.12.8 of \cite{MMO}, $K$ is an ellipsoid. \qed

The version of Theorem \ref{Main} for dimensions $n\geq3$ and any codimension less than $n-1$ is true. The proof follows from our Theorem \ref{Main}, using standard arguments in the literature.

\sn{\bf Acknowledgments.} L. Montejano acknowledges  support  from CONACyT under 
project 166306 and  from PAPIIT-UNAM under project IN112614. E. Morales-Amaya acknowledges  support  from CONACyT, SNI 21120.

\bigskip
\bigskip
\bigskip


\begin{thebibliography}{99}

\bibitem{Bru} Brunn, H., \emph{\"Uber Kurven ohne Wendepunkte.} Habilitationschrift, Ackermann, M\"unchen. 1889.
\bibitem{B1} Burton G.R., Sections of Convex Bodies, J. London Math. Society 12(1976), 331-336
\bibitem{B2} Burton G.R., Some characterizations of the ellipsoid. Israel J. of Math. 28 (1977),339-
\bibitem{L} Larman, D.G., {\it A note on the false center problem}, 
Mathematika 21 (1974), 216-217.

\bibitem{MMO} Martini, H., Montejano, L., Oliveros, D., {\em Bodies of Constant Width; An introduction to convex geometry with applications.} Birkh\"auser, Boston, Bassel, Stuttgart, 2019. 
\bibitem{P} Petty, C. M., Ellipsoids, in: Convexity and its Applications, Eds. P. M. Gruber and J. M. Wills, pp. 264–276, Birkh"auser, Basel, 1983. 
\bibitem{R} Rogers, C.A., {\it Sections and projections of convex bodies,} Portugal Math. 24 (1965), 99-103

\bibitem{So}Soltan, V.,{\em Characteristic properties of ellipsoids and convex quadrics}, Aequat. Math., 93 (2019),371-413



\end{thebibliography}
\end{document}